# ANALYSIS OF TOP TO BOTTOM-$K$ SHUFFLES

By Sharad Goel[1]

*Cornell University*

A deck of $n$ cards is shuffled by repeatedly moving the top card to one of the bottom $k_n$ positions uniformly at random. We give upper and lower bounds on the total variation mixing time for this shuffle as $k_n$ ranges from a constant to $n$. We also consider a symmetric variant of this shuffle in which at each step either the top card is randomly inserted into the bottom $k_n$ positions or a random card from the bottom $k_n$ positions is moved to the top. For this reversible shuffle we derive bounds on the $L^2$ mixing time. Finally, we transfer mixing time estimates for the above shuffles to the lazy top to bottom-$k$ walks that move with probability $1/2$ at each step.

**1. Introduction.** A deck of $n$ cards can be shuffled by repeatedly removing the top card and inserting it uniformly at random back into the deck. A coupling argument shows that the total variation mixing time for this Markov chain is $n \log n$ (see, e.g., [1, 2, 12]). In fact, a detailed analysis yields a closed form expression for the distribution of this chain after any number of steps (see [3]).

Here we analyze a class of walks that generalizes the top to random chain, namely, top to bottom-$k$ shuffles. These shuffles are generated by moving the top card uniformly at random to any of the bottom $k_n$ positions of the deck. For $k_n = n$, we recover the top to random walk. For $k_n = 2$, this is the Rudvalis shuffle, and upper and lower bounds of order $n^3 \log n$ have been shown by Hildebrand [9] and Wilson [15], respectively.

More formally, let $S_n$ be the permutation group, and let $\sigma \in S_n$ denote an element of this group, interpreting $\sigma(i) = j$ to mean that position $i$ holds the card with label $j$. Fix $n \geq k_n > 1$, and denote a cycle permutation by

$$\sigma_l = (1, 2, \ldots, l),$$

Received September 2004; revised June 2005.
[1]Supported in part by NSF Grant DMS-03-06194.
*AMS 2000 subject classifications.* Primary 60; secondary 68.

where $\sigma_l(i) = i+1$ for $1 \leq i \leq l-1$, $\sigma_l(l) = 1$, and $\sigma_l(i) = i$ otherwise. Define the probability measure $q_{n,k_n}$ on $S_n$ by

$$q_{n,k_n}(\sigma) = \begin{cases} \dfrac{1}{k_n}, & \text{if } \sigma = \sigma_l \text{ for some } n - k_n + 1 \leq l \leq n, \\ 0, & \text{otherwise,} \end{cases}$$

and let $\pi$ be the uniform distribution on $S_n$. Then the top to bottom-$k$ shuffle driven by $q_{n,k_n}$ is nonreversible, aperiodic and irreducible with stationary distribution $\pi$.

Let $q^*_{n,k_n}$ denote the bottom-$k$ to top shuffle. It is well known that studying this reversed shuffle is equivalent to studying $q_{n,k_n}$ (see Section 2). Then for the top to bottom-$k$ walk $q_{n,k_n}$, and the reversible variant $\tilde{q}_{n,k_n} = \frac{1}{2}(q_{n,k_n} + q^*_{n,k_n})$, we derive bounds on the total variation and $L^2$ mixing times $T$ and $T_2$. Finally, we show that results for the nonreversible and reversible chains yield as corollaries bound on the lazy top to bottom-$k$ shuffle

$$\hat{q}_{n,k_n} = \tfrac{1}{2}(q_{n,k_n} + \delta_e),$$

where we put weight on the identity.

In particular, our main results are summarized below. In these statements, $A(c)$, $B(c)$, and so on, denote positive, finite constants that may depend on the fixed parameter $c$ but not on $n$.

THEOREM 1.1. *For the top to bottom-$k$ shuffle $q_{n,k_n}$:*

(1) *if $k_n \geq n - \sqrt{(n\log n)/2}$, then*

$$T(S_n, q_{n,k_n}) \sim n\log n;$$

(2) *if $k_n \geq cn$ with $c \in (0,1)$, then*

$$A(c)n\log n \leq T(S_n, q_{n,k_n}) \leq B(c)n^2 \log^2 n;$$

(3) *if $k_n \leq C$, then*

$$A(C)n^3 \leq T(S_n, q_{n,k_n}) \leq B(C)n^3 \log n;$$

(4) *if $k_n = 2, 3$, then*

$$An^3 \log n \leq T(S_n, q_{n,k_n}) \leq Bn^3 \log n.$$

For Theorem 1.1(1), the notation $\sim$ indicates that the walk presents a total variation cut-off at time $n \log n$. See Lemma 3.2 for a precise statement of the result. After this paper was submitted, the author learned of the work of Jonasson [10], who shows that nonreversible top to bottom-$k$ shuffles have total variation mixing time $T(S_n, q_{n,k_n}) \approx n^3 \log n / k_n^2$ uniformly for all choices of $k_n$.

THEOREM 1.2. *Let $\tilde{q}_{n,k_n} = \frac{1}{2}(q_{n,k_n} + q^*_{n,k_n})$ be the additive symmetrization of the top to bottom-k shuffle. Then:*

(1) *if $k_n \geq n - C$, then*
$$T(S_n, \tilde{q}_{n,k_n}) \leq T_2(S_n, \tilde{q}_{n,k_n}) \leq B(C) n \log n;$$

(2) *if $k_n \leq cn$ with $c \in (0,1)$, then*
$$T_2(S_n, \tilde{q}_{n,k_n}) \geq T(S_n, \tilde{q}_{n,k_n}) \geq A(c) n^2;$$

*and*
$$T_2(S_n, \tilde{q}_{n,k_n}) \geq \frac{A(c) n^3}{k_n^2} \log n;$$

(3) *for any $k_n$,*
$$T(S_n, \tilde{q}_{n,k_n}) \leq T_2(S_n, \tilde{q}_{n,k_n}) \leq B n^3 \log n.$$

*In particular, (2) and (3) show that if $k_n \leq C$, then*
$$A(C) n^3 \log n \leq T_2(S_n, \tilde{q}_{n,k_n}) \leq B n^3 \log n.$$

The two lower bounds in Theorem 1.2(2) are complimentary in the sense that the first gives better estimates for $k_n \approx cn$, while the second works best for $k_n \ll cn$.

THEOREM 1.3. *For the lazy top to bottom-k shuffle $\hat{q}_{n,k_n}$:*

(1) *if $k_n \geq n - C$, then*
$$An \log n \leq T_2(S_n, \hat{q}_{n,k_n}) \leq B(C) n \log n;$$

(2) *if $k_n \geq n - \sqrt{(n \log n)/2}$, then*
$$T(S_n, \hat{q}_{n,k_n}) \sim 2n \log n;$$

(3) *if $k_n \geq cn$ with $c \in (0,1)$, then*
$$A(c) n \log n \leq T(S_n, \hat{q}_{n,k_n}) \leq B(c) n^2 \log^2 n;$$

(4) *if $k_n = 2, 3$, then*
$$An^3 \log n \leq T(S_n, \hat{q}_{n,k_n}) \leq B n^3 \log n;$$

(5) *for any $k_n$,*
$$T_2(S_n, \hat{q}_{n,k_n}) \leq B n^3 \log n.$$

For Theorem 1.3(2), the notation $\sim$ indicates that the walk presents a total variation cut-off at time $2n \log n$. See Remark 5.1 for a precise statement of the result. Also observe that the estimates of Theorem 1.3(1)(2) bound the $L^2$ mixing time $T_2$ and the total variation mixing time $T$, respectively.

As $k_n$ varies from a constant to $n$, these results are most satisfactory at the extremes of the range. For large $k_n$ the walks behave like the top to random chain, mixing in $n \log n$ steps. Theorem 1.1(1) proves mixing in the strongest possible sense: cut-off at precisely $n \log n$. Let us note here that the precise $L^2$ cut-off time is not yet known even for the top to random shuffle $q_{n,n}$.

For small $k_n$, the walks behave like the Rudvalis shuffle, mixing in $n^3 \log n$ steps. Theorem 1.2 proves this for the reversible chain, whereas Theorems 1.1 and 1.3 give complete results only for $k_n = 2, 3$.

The worst gap in these results occurs when $k_n \approx n/2$. For these "top to bottom half" shuffles, [10] shows a $\Theta(n \log n)$ mixing time for the nonreversible shuffle, and our results give an $\Omega(n^2)$ lower bound for the reversible shuffle. In particular, the nonreversible and reversible top to bottom half shuffles mix at different rates. In this range, one difficulty in analyzing the reversible walk is that comparison with random transposition, one of the best understood models of random walk, can at best yield $O(n^3 \log n)$ upper bounds (see Lemma 4.5).

A variety of methods are used to prove the results of this paper. The upper bounds for the nonreversible top to bottom-$k$ shuffle are found by coupling arguments. The lower bound in Theorem 1.1(4) uses Wilson's lemma (see, e.g., [13, 15]). For the reversible chain, we use comparison techniques for walks on finite groups to prove both upper and lower bounds (see, e.g., [4]). Notably, comparison previously has been applied only to find upper bounds. It appears that this is the first application of comparison techniques to prove lower bounds.

In Section 2 we introduce our notation and review basic Markov chain theory. Sections 3 and 4 give proofs of Theorems 1.1 and 1.2, respectively. Finally, Section 5 applies results from the previous sections to find bounds on the lazy walk $\hat{q}_{n,k_n}$.

**2. Notation and basics.** Let $G$ be a finite group with probability measure $q$, and let $\{\eta_i\}$ be $G$-valued independent random variables with distribution $q$. The left-invariant walk on $G$ driven by $q$ is defined by $X_0 = e$ and

$$X_{k+1} = X_k \cdot \eta_k.$$

For the top to bottom-$k$ measure $q_{n,k_n}$, this definition corresponds to the informal card shuffling description given in the Introduction. See, for example, [12] for more details. Define convolution powers of $q$ by

$$q^m(g) = q^{m-1} \star q(g) = \sum_{h \in G} q^{m-1}(h) q(h^{-1} g).$$

Equivalently, $q^m(g)$ can be thought of as the sum of weighted paths:

$$q^m(g) = \sum \prod_{i=1}^{m} q(h_i),$$

where, for each fixed $g$, the sum is taken over $m$-tuples $(h_1, \ldots, h_m)$ such that $h_1 \cdots h_m = g$.

If $\mathrm{supp}(q) = \{g : q(g) > 0\}$ is not contained in a proper subgroup of $G$ or in a coset of a proper normal subgroup, then

$$q^m(g) \to \frac{1}{|G|} \quad \text{as } m \to \infty.$$

Our results give bounds on the rate of convergence for $q_{n,k_n}$ and its variants. Intuitively, these results are estimates on the number of top to bottom-$k$ shuffles needed to mix a deck of cards. To make this statement more precise, we first need a way to measure distance between the distribution of the chain at time $m$ and the stationary distribution. For probability measures $q$ and $\pi$ on a finite group $G$, define the total variation distance as

$$\|q - \pi\|_{\mathrm{TV}} = \sup_{A \subset G} |q(A) - \pi(A)| = \tfrac{1}{2} \sum_{g \in G} |q(g) - \pi(g)|.$$

Alternatively, some of our results will be in terms of the $L^p(\pi)$ distance

$$d_{\pi,p}(q) = \left\| \frac{q}{\pi} - 1 \right\|_{L^p(\pi)} = \left( \sum_{g \in G} \left| \frac{q(g)}{\pi(g)} - 1 \right|^p \pi(g) \right)^{1/p}.$$

Observe that $\|q - \pi\|_{\mathrm{TV}} = \tfrac{1}{2} d_{\pi,1}(q)$. Our results are for the cases $p = 1, 2$.

We define the deck to be shuffled when the distance between the distribution of the deck and the stationary distribution is small. Namely, the total variation mixing time is given by

$$T(G, q) = \inf\left\{ m \,\big|\, \|q^m - \pi\|_{\mathrm{TV}} \leq \frac{1}{2e} \right\}$$

and the $L^p$ mixing time by

$$T_p(G, q) = \inf\left\{ m \,\Big|\, \left\| \frac{q^m}{\pi} - 1 \right\|_{L^p} \leq \frac{1}{e} \right\}.$$

With these definitions, by Jensen's inequality,

$$T(G, q) = T_1(G, q) \leq T_2(G, q).$$

Moreover, the functions $k \mapsto d_{\pi,p}(q^k)$ are nonincreasing and sub-additive. In particular, for $k \geq T_p(G, q)$,

$$d_{\pi,p}(q^k) \leq e^{-\lfloor k/T_p(G,q) \rfloor}.$$

This inequality motivates our somewhat arbitrary choice of $1/e$ in the definition of mixing time. For details, see, for example, [6, 12].

The Markov operator $Q$ associated to a probability measure $q$ on $G$ is given by $Qf = f \star q^*$, where $q^*(g) = q(g^{-1})$. The reversed random walk is driven by $q^*$ and has as its associated operator the adjoint of $Q$. That is, $q^*$ has associated Markov operator $Q^*f = f \star q$.

Note that since we are on a group, the stationary measure $\pi$ is uniform, and furthermore,

$$\sum_{g \in G} \left| \frac{q(g)}{\pi(g)} - 1 \right|^p \pi(g) = \sum_{g \in G} \left| \frac{q(g^{-1})}{\pi(g)} - 1 \right|^p \pi(g).$$

Consequently, $d_{\pi,p}(q) = d_{\pi,p}(q^*)$, and with respect to analyzing mixing time, we can study either the walk or its reversal.

For a sequence of numbers $a_n, b_n$, we use the notation $a_n \preceq b_n$ to indicate that there is a universal constant $C > 0$ (independent of $n$) such that $a_n \leq Cb_n$. For the two-sided bound, we use $a_n \approx b_n$ to indicate that there are constants $c, C > 0$ such that $ca_n \leq b_n \leq Ca_n$. For mixing times $T(n)$, the notation $T(n) \sim a_n$ indicates cut-off at time $a_n$. For a precise definition of cut-off, see, for example, [11].

**3. Proof of Theorem 1.1.** In this section we present upper and lower bounds for the mixing time of the nonreversible walk $q_{n,k_n}$, using primarily probabilistic techniques. For $k_n = 3$, we use the method of [15] to derive a lower bound.

To prove mixing time bounds for the top to bottom-$k$ shuffle, we make extensive use of the following well-known coupling result (see, e.g., [1, 2, 12]).

THEOREM 3.1. *Let $q$ be a probability measure on a finite group $G$. Let $(X_n^1, X_n^2)$ be a coupling for the random walk driven by $q$ with $(X_n^1)$ starting at the identity and $(X_n^2)$ starting from the stationary distribution $\pi$ [i.e., $\mathrm{dist}(X_0^2) = \pi$]. Then*

$$\|q^m - \pi\|_{\mathrm{TV}} \leq \mathbb{P}(T > m),$$

*where*

$$T = \inf\{m | \forall\, k \geq m,\ X_k^1 = X_k^2\}.$$

*Furthermore, there exists a coupling such that the inequality above is an equality.*

We will also make use of the following coupon-collectors lemma (see, e.g., [1]).

LEMMA 3.1 (Coupon-collectors lemma). *Let $R_m$ be the number of distinct cards obtained in $m$ uniform random draws with replacement from a deck of $n$ cards. That is, $R_m = |\{C_1, \ldots, C_m\}|$ with $C_i$ i.i.d. uniform on $\{1, \ldots, n\}$. Let $L_j = \min\{m | R_m = n - j\}$, that is, the number of draws before all but $j$ cards have been chosen. Then for fixed $j$,*

$$\frac{L_j}{n \log n} \to 1 \quad \text{in probability.}$$

In the case of $q_{n,n}$, that is, the top to random shuffle, the correct mixing time $n \log n$ can be found using a coupling of the time reversed process $q_{n,n}^*$. For this random to top shuffle, the coupling is as follows: choose a label uniformly at random from $\{1, \ldots, n\}$ and in each deck move the card with this label to the top. Clearly, this is a coupling, and the coupling time is given by the coupon-collectors lemma (for details, see, e.g., [1]). The proof of Lemma 3.2 is by a similar coupling.

LEMMA 3.2. *For $k_n \geq n - \sqrt{\frac{1}{2} n \log n}$, the walk $(S_n, q_{n,k_n})$ presents a total variation cut-off at $t_n = n \log n$. That is, for $\varepsilon \in (0, 1)$,*

$$\lim_{n \to \infty} \|q_{n,k_n}^{(1+\varepsilon)n \log n} - \pi\|_{\text{TV}} = 0$$

*and*

$$\lim_{n \to \infty} \|q_{n,k_n}^{(1-\varepsilon)n \log n} - \pi\|_{\text{TV}} = 1.$$

PROOF. Since $d_{\text{TV}}(p^{(n)}, u) = d_{\text{TV}}(p^{*(n)}, u)$, we can consider the reversed random walk $q_{n,k_n}^*$. For this reversed walk, we define a coupling $(X_1^m, X_2^m)$, where $X_1$ starts from the identity and $X_2$ is drawn from the stationary distribution. Let

$$A_j^m = \{X_j^m(i) | n - k_n < i \leq n\}, \qquad j = 1, 2.$$

That is, $A_j^m$ is the set of cards that at time $m$ are in the bottom $k_n$ positions of deck $j$. At time $m$, in the first deck pick a card $\sigma_a$ uniformly at random from $A_1^m$ and move it to the top of the deck. If $\sigma_a \in A_2^m$, then move card $\sigma_a$ in the second deck to the top. If not, then in the second deck uniformly at random pick a card from $A_2^m \setminus A_1^m$ and move it to the top.

Clearly, deck one is driven by $q_{n,k_n}^*$. For the second deck, note that any card in $A_1^m \cap A_2^m$ is chosen if and only if it is chosen in the first deck, and hence, with probability $1/k_n$. And cards in $A_2^m \setminus A_1^m$ are chosen with probability

$$\frac{k_n - |A_1^m \cap A_2^m|}{k_n} \cdot \frac{1}{k_n - |A_1^m \cap A_2^m|} = \frac{1}{k_n}.$$

So this is, in fact, a coupling. Define

$$\tau_0 = \inf\{m | X_1^m(i) = X_2^m(i) \text{ for } 1 \leq i \leq n - k_n\}.$$

That is, $\tau_0$ is the first time the top $n - k_n$ cards are matched in both decks. Then for $m > \tau_0$, $A_1^m = A_2^m$, that is, the set of cards in the bottom $k_n$ positions are the same in each deck. Consequently, after time $\tau_0$, new matches are not broken and every time an unmatched card is chosen, a new match is made.

First we estimate $\tau_0$. Let $L$ be the probability that, starting with all cards unmatched, $n - k_n$ consecutive matches are made. Then,

$$L \geq \left(1 - \frac{n - k_n}{k_n}\right)^{n - k_n}$$

$$\geq \left(1 - \frac{1}{\sqrt{n}/\sqrt{(1/2)\log n} - 1}\right)^{\sqrt{(1/2)n \log n}}$$

$$\approx \frac{1}{\sqrt{n}}.$$

Furthermore, by the Markov property, for fixed $\varepsilon > 0$,

$$P(\tau_0 \geq \varepsilon n \log n) \leq P\left(\tau_0 \geq \varepsilon \sqrt{n \log n} \cdot \sqrt{\frac{1}{2} n \log n}\right)$$

$$\leq \left[1 - \left(1 - \frac{1}{\sqrt{n}/\sqrt{(1/2)\log n} - 1}\right)^{\sqrt{(1/2)n \log n}}\right]^{\varepsilon \sqrt{n \log n}}$$

$$\xrightarrow{n \to \infty} 0.$$

Let $\tau_1$ be the time it takes after $\tau_0$ for each card in $A_1^{\tau_0} = A_2^{\tau_0}$ to be selected. That is,

$$\tau_1 = \inf\{m | m > 0, \text{ each card in } A_1^{\tau_0} \text{ has been selected by time } m + \tau_0\}.$$

By the coupon-collectors lemma, for fixed $\varepsilon > 0$,

$$\lim_{n \to \infty} P(\tau_1 \geq (1 + \varepsilon) k_n \log k_n) = 0.$$

Finally, if $T$ is the coupling time, then since

$$P(T > (1 + \varepsilon) n \log n) \leq P\left(\tau_0 \geq \frac{\varepsilon}{2} n \log n\right) + P\left(\tau_1 \geq \left(1 + \frac{\varepsilon}{2}\right) n \log n\right)$$

$$\xrightarrow{n \to \infty} 0$$

by Theorem 3.1,

$$\lim_{n \to \infty} \|q_{n,k_n}^{(1+\varepsilon)n \log n} - \pi\|_{\text{TV}} = 0.$$

The lower bound argument is analogous to that of the top to random shuffle (see, e.g., [1]). Let $B_j$ be the set of permutations $\sigma$ for which the bottom $j$ cards have increasing labels. That is,
$$\sigma(n-j+1) < \sigma(n-j+2) < \cdots < \sigma(n).$$
Then $\pi(B_j) = \frac{1}{j!}$. Starting from the identity, let $L_j$ be the number of shuffles until all but $j$ of the cards with labels in $\{n-k_n+1, \ldots, n\}$ have been chosen. Then, if $L_j > m$, the bottom $j$ cards after $m$ bottom $k_n$ to top shuffles are in increasing order. So for fixed $\varepsilon > 0$, there is an $\varepsilon' > 0$ such that
$$\lim_{n\to\infty} \|q_{n,k_n}^{(1-\varepsilon)n\log n} - \pi\|_{\mathrm{TV}} \geq \lim_{n\to\infty} P(L_j > (1-\varepsilon)n\log n) - \frac{1}{j!}$$
$$\geq \lim_{n\to\infty} P(L_j > (1-\varepsilon')k_n \log k_n) - \frac{1}{j!}$$
since
$$\lim_{n\to\infty} \frac{k_n \log k_n}{n \log n} = 1.$$
Using the coupon-collectors lemma, the result follows. $\square$

Lemmas 3.3 and 3.4 below bound the mixing time of $q_{n,k_n}$ in the cases where $k_n$ is relatively large and when $k_n$ is small. Both lemmas rely on the following coupling.

We construct a coupling $(X_1^m, X_2^m)$ where $X_1$ starts from the identity and $X_2$ is drawn from the uniform distribution. Recall that the notation $X_s^m(i) = j$ can be interpreted to mean that at time $m$ position $i$ in deck $s$ holds the card with label $j$. Let
$$A_s^m = \{X_s^m(i) | n - k_n + 2 \leq i \leq n\}, \qquad s = 1, 2.$$
Note that $A_s^m$ is not the set of cards in the bottom $k_n$ positions (to which the top card can be sent), but rather only the cards in the bottom $k_n - 1$ positions.

We define a coupling as follows: first pick one of the two decks with equal probability. Say we picked deck one. Then $X_1$ proceeds as usual by uniformly at random moving the top card to one of the bottom $k_n$ positions; $X_2$ mimics the moves of $X_1$ except in a couple of cases. If $X_1^m(1) \in A_2^m$ (i.e., the top card in the first deck is in $A_2^m$), and the first deck moves the top card to position $(X_2^m)^{-1}(X_1^m(1))$, then the second deck moves the top card to $(X_2^m)^{-1}(X_1^m(1)) - 1$. And, if $X_1^m(1) \in A_2^m$ and the first deck moves the top card to $(X_2^m)^{-1}(X_1^m(1)) - 1$, then the second deck moves the top card to $(X_2^m)^{-1}(X_1^m(1))$. We have an analogous description if we originally picked deck two. Accordingly, if card $i$ is on the top of one deck and in the bottom $k_n - 1$ positions of the other deck, then, with probability $1/k_n$, it will couple on the next move. Furthermore, matches between the decks are never broken.

LEMMA 3.3. *For $c \in (0,1)$ and $k_n \geq cn$, there exist constants $A(c)$ such that the total variation mixing time for the walk driven by $q_{n,k_n}$ satisfies*

$$T(S_n, q_{n,k_n}) \leq An^2 \log^2 n.$$

PROOF. We use the coupling described above. Let $\tau_j$ be the first time that the cards with label $j$ couple in the two decks. That is,

$$\tau_j = \inf\{m | (X_1^m)^{-1}(j) = (X_2^m)^{-1}(j)\}.$$

We estimate $\tau_j$ by first showing that, starting from any permutation of the decks, any card $j$ has probability at least $\frac{C}{n}$ to couple within $3n \log n$ steps. Let $\tau_\sigma^j$ be the first time card $j$ reaches the top of deck one, starting from state $\sigma$. And let $\tau_\sigma$ be the first time the bottom card reaches position $n - k_n$. Then for $n$ sufficiently large,

$$P(\tau_\sigma^j > 2n \log n) \leq P(\tau_\sigma > 2n \log n - (n - k_n))$$

$$\leq k_n \exp\left(-\frac{2n \log n - (n - k_n)}{k_n}\right)$$

$$\leq \frac{1}{2}.$$

The second inequality follows from the fact that $\tau_\sigma$ is the sum of independent geometric waiting times with means $k_n, k_n/2, \ldots, k_n/k_n$, and, consequently, is equivalent to the coupon collectors problem. In particular, the above shows that, starting from any state, there is positive probability independent of $n$ and $k_n$ that card $j$ reaches the top of the first deck in $2n \log n$ steps.

When card $j$ gets to the top of the first deck, we are in one of three situations: card $j$ is already coupled, card $j$ in the second deck is in the bottom $k_n - 1$ positions, or card $j$ in the second deck is in the top $n - k_n + 1$ positions. In the first two situations, card $j$ will be coupled at the next step with probability at least $1/k_n$ (if $j$ is already coupled, it will remain coupled at the next step). So we only need to consider the third situation. Assume card $j$ moves to one of the bottom $\lceil Bk_n \rceil$ positions for some $B \in (0,1)$ (which happens with probability at least $B$). Let $\tau_B$ be the first time $j$ leaves the bottom $k_n - 1$ positions. Then $\tau_B$ is the sum of independent geometric waiting times, and depends on the exact position in the bottom $\lceil Bk_n \rceil$ to which card $j$ moves. However, by construction, we have the lower bound

$$E\tau_B \geq \sum_{r=\lceil Bk_n \rceil}^{k_n-1} \frac{k_n}{r}$$

$$\geq k_n \log \frac{1}{B + 1/k_n}.$$

And,
$$\operatorname{Var}(\tau_B) \leq \sum_{r=1}^{k_n-1} \frac{k_n(k_n-r)}{r^2}$$
$$\leq 2k_n^2.$$

By Chebyshev's inequality,
$$P\left(\tau_B > \frac{E\tau_B}{2}\right) \geq 1 - \frac{4\operatorname{Var}(\tau_B)}{(E\tau_B)^2}$$
$$\geq 1 - \frac{8}{\log^2 1/(B+1/k_n)}.$$

Consequently, if we choose $B$ and $K$ such that
$$\log \frac{1}{B+1/K} \geq \max\left(\frac{2(1-c)}{c}, 3\right),$$
where $c$ is from the statement of the lemma, then there exists $\delta > 0$ (independent of $n$) such that, for $k_n \geq K$,
$$P(\tau_B > n - k_n) \geq P(\tau_B > E\tau_B/2) \geq \delta.$$

For instance, we can choose $\delta = 1/9$. But if $\tau_B > n - k_n$, then $j$ will still be in the bottom $k_n - 1$ positions of deck one when $j$ reaches the top of deck two. Consequently, for each of the original three cases, after reaching the top of deck one, card $j$ couples within the next $n - k_n$ steps with probability at least $\delta/n$. Combining this with the bound on $\tau_\sigma^j$, for the coupling time $\tau_j$ of card $j$, we have
$$P(\tau_j \leq 3n \log n) \geq \frac{\delta}{2n}.$$

Moreover, by the Markov property,
$$P(\tau_j > An^2 \log^2 n) \leq \left(1 - \frac{\delta}{2n}\right)^{An \log n/3}$$
$$\leq \exp\left(-\frac{\delta A \log n}{6}\right).$$

Finally, if $T$ is the coupling time for the two decks, then
$$P(T > An^2 \log^2 n) \leq n \exp\left(-\frac{\delta A \log n}{6}\right)$$

and the result follows by taking $A$ sufficiently large. $\square$

REMARK 3.1. Using the lower bound argument of Lemma 3.2, we can show that, for $c \in (0,1)$, $k_n \geq cn$, there exist constants $B(c)$ such that the mixing time satisfies

$$T(S_n, q_{n,k_n}) \geq B(c) n \log n.$$

The following lemma gives an upper bound on the mixing time for the walk driven by $q_{n,k_n}$ with $k_n \leq C$. The coupling used to prove the result is the same as in Lemma 3.3, however, we analyze the coupling time by a different technique.

LEMMA 3.4. *For $k_n \leq C$, there exist constants $A(C)$ such that the total variation mixing time for the walk driven by $q_{n,k_n}$ satisfies*

$$T(S_n, q_{n,k_n}) \leq A n^3 \log n.$$

PROOF. Using the coupling described above, we show that, starting from any permutation of the decks, any card $i$ has probability at least $\delta > 0$ (independent of $n$) to couple within $n^3$ steps. Fix card $i$ and let $\tau$ be the first time that card $i$ is on the top of one deck and in the bottom $k_n - 1$ positions of the other. Then at the next step, the cards have probability $1/k_n$ to couple. Let

$$\tau_j^1 = \inf\{t | X_j^t(1) = i\},$$
$$\tau_j^m = \inf\{t > \tau_j^{m-1} | X_j^t(1) = i\}.$$

That is, $\tau_j^m$ is the time when card $i$ is on top of deck $j$ for the $m$th time. Without loss of generality, assume that $\tau_1^1 \leq \tau_2^1$. If $\tau_j^m \leq \tau$, then

$$\tau_1^m \leq \tau_2^m \leq \tau_1^m + n - k_n.$$

And if $\tau_j^{m+1} \leq \tau$, then

$$\tau_j^{m+1} \leq \tau_j^m + 2(n - k_n).$$

Define the random variables $d_i^m = [(X_1^m)^{-1}(i) - (X_2^m)^{-1}(i)] \mod n$, which give the oriented distance between the positions of the $i$th card in each deck. Note that $d_i^m$ only changes when $i$ is in the bottom $k_n - 1$ positions in at least one deck. Let $\tau^* = \inf\{t > \tau_1^1 | X_1^t(i) = n - k_n + 1\}$. Then define $Y_h^l$ as i.i.d. random variables with distribution given by

$$P(Y_h^l = t) \stackrel{\text{def}}{=} P(\tau^* - \tau_1^1 = t).$$

That is, $Y_h^l$ gives the amount of time it takes a card to get from the top of the deck to the $n - k_n + 1$ position. Furthermore, before $\tau$, the distribution of the change in distance is given by $d_i^{\tau_1^{m+1}} - d_i^{\tau_1^m} \stackrel{\text{dist}}{=} Y_1^m - Y_2^m$. Consequently,

$$P(\tau_1^{m+1} \leq \tau) \leq P\left(\left|\sum_{l=1}^m Y_1^l - Y_2^l\right| \leq n\right).$$

Let $\sigma^2 = \text{Var}(Y_1^l - Y_2^l)$, and note that $\sigma < \infty$ since $Y_h^l$ can be realized as a finite sum of geometric waiting times. Then by the central limit theorem, by taking $m = n^2$, we have that $P(\tau_1^{n^2+1} \leq \tau) \leq 1 - \varepsilon$ independent of $n$. That is, $P(\tau < \tau_1^{n^2+1}) > \varepsilon$. Furthermore, since $\tau_1^1 \leq n$ with positive probability independent of $n$, $P(\tau < 3n^3) > \varepsilon$. Consequently, there is a $\delta > 0$ such that if $\tau_i$ is the coupling time for card $i$, then $P(\tau_i < 3n^3) > \delta$. Finally, if $T$ is the coupling time for the two decks, then

$$\begin{aligned} P(T > An^3 \log n) &\leq nP(\tau_i > An^3 \log n) \\ &\leq n(1-\delta)^{A \log n/3} \\ &\stackrel{A \to \infty}{\longrightarrow} 0. \end{aligned}$$

By taking $A$ sufficiently large, the result follows from Theorem 3.1. □

REMARK 3.2. For $k_n \leq C$, the walk performed by one card under the measure $q_{n,k_n}$ is an example of a class of walks known as necklace chains. By results in [14], this immediately yields the lower bound

$$B(C)n^3 \leq T(S_n, q_{n,k_n}).$$

In [9], the Rudvalis shuffle $q_{n,2}$ is shown to have an upper bound of order $O(n^3 \log n)$. In [15], a matching lower bound for this shuffle is given by using Theorem 3.2. Here we show that the method of [15] can also used to lower bound $q_{n,3}$.

Given a chain $X_t$, we say that the chain $(\tilde{X}_t, Y_t)$ is a lifting of the original chain if the marginal distribution of $\tilde{X}_t$ is the same as the distribution of $X_t$.

THEOREM 3.2. *Suppose that a Markov chains $X_t$ has a lifting $(X_t, Y_t)$, and that $\Psi$ is an eigenfunction of the lifted Markov chain: $E[\Psi(X_{t+1}, Y_{t+1})|(X_t, Y_t)] = \lambda \Psi(X_t, Y_t)$. Suppose that $|\Psi(x,y)|$ is a function of $x$ alone, $|\lambda| < 1$, $\Re(\lambda) \geq 1/2$, and that we have an upper bound $R$ on $E[|\Psi(X_{t+1}, Y_{t+1}) - \Psi(X_t, Y_t)|^2|(X_t, Y_t)]$. Let $\gamma = 1 - \Re(\lambda)$. Then when the number of steps $t$ is bound by*

$$t \leq \frac{\log \Psi_{\max} + (1/2) \log \gamma \varepsilon/(4R)}{-\log(1-\gamma)},$$

*the variation distance satisfies $\|X_t - \pi\|_{\text{TV}} \geq 1 - \varepsilon$.*

For a discussion of Theorem 3.2, see [13, 15, 16].

LEMMA 3.5. *For $\varepsilon > 0$, there exist constants $C(\varepsilon), N > 0$ such that, for $n \geq N$,*

$$\|q_{n,3}^m - \pi\|_{\mathrm{TV}} \geq 1 - \varepsilon$$

*for $m \leq Cn^3 \log n$.*

PROOF. Let $X_t^{-1}(j) = j'$ denote that the card with label $j$ is at position $j'$ at time $t$. First we lift the chain to $(X_t^{-1}, Y_t) = (X_t^{-1}, t \bmod n)$. Let $Z_t(j) = (X_t^{-1}(j) - X_0^{-1}(j) + Y_t(j)) \bmod n$ and let $\eta(t) \in \{\sigma_{n-2}, \sigma_{n-1}, \sigma_n\}$ denote the cycle that is chosen at time $t$. Then,

$$(X_{t+1}^{-1}(j), Z_{t+1}(j)) = \begin{cases} (X_t^{-1}(j), Z_t(j) + 1), & \eta_t = \sigma_{n-1}, X_t^{-1}(j) = n \text{ or} \\ & \eta_t = \sigma_{n-2}, X_t^{-1}(j) \geq n-1, \\ (X_t^{-1}(j) - 1, Z_t(j)), & \eta_t = \sigma_n \text{ or} \\ & \eta_t = \sigma_{n-1}, X_t^{-1}(j) \leq n-1 \text{ or} \\ & \eta_t = \sigma_{n-2}, X_t^{-1}(j) \leq n-2, \\ (n-1, Z_t(j) - 1), & \eta_t = \sigma_{n-1}, X_t^{-1}(j) = 1, \\ (n-2, Z_t(j) - 2), & \eta_t = \sigma_{n-2}, X_t^{-1}(j) = 1. \end{cases}$$

Define $v(x)$ to be the $x$th number in the list

$$\lambda^{n-3}, \ldots, \lambda, 1, \chi_1, \chi_0$$

and define the functions

$$\Psi_j(X_t^{-1}, Y_t) = v(X_t^{-1}(j)) \exp(Z_t(j) 2\pi i/n),$$

$$\Psi(X_t^{-1}, Y_t) = \sum_{j=1}^n \Psi_j(X_t^{-1}, Y_t).$$

Now we will find values for $\lambda, \chi_1, \chi_0$ that make $\Psi_j$ (and, hence, $\Psi$) an eigenfunction. Also note that $|\Psi(X_t^{-1}, y_1)| = |\Psi(X_t^{-1}, y_2)|$ for all $y_1, y_2$. If $2 \geq X_t^{-1}(j) \geq n-2$, then

$$\Psi_j(X_{t+1}^{-1}, Y_{t+1}) = \lambda \Psi_i(X_t^{-1}, Y_t).$$

Let $w = e^{2\pi i/n}$. By looking at what happens when $X_t^{-1}(j) = 1$, $X_t^{-1}(j) = n$, and $X_t^{-1}(j) = n-1$, we find that $\Psi_j$ is an eigenfunction with eigenvalue $\lambda$ when the equations

$$\chi_0 + \chi_1 w^{-1} + w^{-2} = 3\lambda^{n-2},$$

$$\frac{\chi_1}{\chi_0} + 2w = 3\lambda,$$

$$\frac{2}{\chi_1} + w = 3\lambda$$

are satisfied. In particular,
$$\chi_0 = \frac{2}{(3\lambda - w)(3\lambda - 2w)},$$
$$\chi_1 = \frac{2}{3\lambda - w}$$

and $\lambda$ is a root of the polynomial
$$f(\lambda) = 9\lambda^n - 9w\lambda^{n-1} + 2w\lambda^{n-2} - 3w^{-2}\lambda^2 + w^{-1}\lambda.$$

We will use Newton's method to approximate a root of $f(\lambda)$ starting with $z_0 = 1$ and $z_{k+1} = z_k - f(z_k)/f'(z_k)$. By Taylor's theorem,
$$|f(z_{k+1})| \leq \frac{1}{2} \max_{0 \leq p \leq 1} |f''(pz_k + (1-p)z_{k+1})| \cdot \left|\frac{f(z_k)}{f'(z_k)}\right|^2.$$

Furthermore, since
$$f'(\lambda) = 9n\lambda^{n-1} - 9(n-1)w\lambda^{n-2} + 2w(n-2)\lambda^{n-3} - 6w^{-2}\lambda + w^{-1},$$
$$f''(\lambda) = 9n(n-1)\lambda^{n-2} - 9(n-1)(n-2)w\lambda^{n-3}$$
$$+ 2w(n-2)(n-3)\lambda^{n-4} - 6w^{-2}$$

if $z = 1 + O(1/n^2)$, then $f'(z) = 2n + O(1)$ and $f''(z) = 2n + O(n)$. So if $z_k = 1 + O(1/n^2)$ and $z_{k+1} = 1 + O(1/n^2)$, then
$$|f(z_{k+1})| \leq \frac{1 + O(1/n)}{4} |f(z_k)|^2.$$

Furthermore,
$$f(z_0) = 9 - 7w + w^{-1} - 3w^{-2}$$
$$= \frac{36\pi^2}{n^2} - i\frac{4\pi}{n} + O(1/n^4).$$

Consequently, by induction,
$$|f(z_k)| \leq 4\left(\frac{\pi}{n}\right)^{2^k} + O\left(\frac{1}{n^{2^k+1}}\right)$$
$$|z_{k+1} - z_k| = \frac{2}{n}\left(\frac{\pi}{n}\right)^{2^k} + O\left(\frac{1}{n^{2^k+2}}\right).$$

So for $n$ sufficiently large, the sequence $\{z_k\}$ converges to a point $z_\infty$ and by continuity, $f(z_\infty) = 0$. Furthermore, since
$$f'(z_0) = 9n - 9(n-1)w + 2w(n-2) - 6w^{-2} + w^{-1}$$
$$= 2n - i14\pi + O(1/n),$$

$$\operatorname{Re}(z_1) = 1 - \operatorname{Re}\left(\frac{f(z_0)}{f'(z_0)}\right)$$

$$= 1 - \frac{\operatorname{Re}(f(z_0))\operatorname{Re}(f'(z_0)) + \operatorname{Im}(f(z_0))\operatorname{Im}(f'(z_0))}{|f'(z_0)|^2}$$

$$= 1 - \left(\frac{18\pi^2 + 14\pi}{n^3}\right) + O(1/n^4).$$

Finally, since

$$|z_1 - z_\infty| \leq \frac{2\pi^2}{n^3} + O(1/n^4),$$

there exist $c_2 > c_1 > 0$ such that

$$1 - \frac{c_1}{n^3} + O(1/n^4) \geq \operatorname{Re}(z_\infty) \geq 1 - \frac{c_2}{n^3} + O(1/n^4).$$

With $\lambda = z_\infty$, $\chi_0 = 1 + O(1/n)$, and $\chi_1 = 1 + O(1/n)$. Consequently,

$$\Psi_{\max} = n + O(1/n).$$

Now we estimate $R$. Since $|\lambda - 1| = O(1/n^2)$,

$$\frac{\Psi_i(X_{t+1}^{-1}, Y_{t+1}) - \Psi_i(X_t^{-1}, Y_t)}{w^{Z_t(i)}}$$

$$= \begin{cases} (\lambda - 1)\lambda^{X_t^{-1}(i)} = O(1/n^2), & 2 \leq X_t^{-1}(i) \leq n-2, \\ \chi_0 - \lambda^{n-3} = O(1/n), & X_t(i) = 1, \eta_t = \sigma_n, \\ \chi_1 w^{-1} - \lambda^{n-3} = O(1/n), & X_t(i) = 1, \eta_t = \sigma_{n-1}, \\ w^{-2} - \lambda^{n-3} = O(1/n), & X_t(i) = 1, \eta_t = \sigma_{n-2}, \\ \chi_1 - \chi_0 = O(1/n), & X_t(i) = n, \eta_t = \sigma_n, \\ w\chi_0 - \chi_0 = O(1/n), & X_t(i) = n, \eta_t = \sigma_{n-1}, \\ w\chi_0 - \chi_0 = O(1/n), & X_t(i) = n, \eta_t = \sigma_{n-2}, \\ 1 - \chi_1 = O(1/n), & X_t(i) = n-1, \eta_t = \sigma_n, \\ 1 - \chi_1 = O(1/n), & X_t(i) = n-1, \eta_t = \sigma_{n-1}, \\ w\chi_1 - \chi_1 = O(1/n), & X_t(i) = n-1, \eta_t = \sigma_{n-2}. \end{cases}$$

Consequently,

$$|\Psi(X_{t+1}^{-1}, Y_{t+1}) - \Psi(X_t^{-1}, Y_t)| = O(1/n)$$

and we can take $R = O(1/n^2)$. The result follows by Theorem 3.2. □

**4. Proof of Theorem 1.2.** In this section we focus on the reversible walk $\frac{1}{2}(q_{n,k_n} + q^*_{n,k_n})$. For reversible chains, path comparison is a useful technique for studying rates of convergence (see, e.g., [4, 5, 6, 8]). In particular, many of the arguments in this section rely on the notion of a flow to compare top to bottom-$k$ shuffles with the well-studied random transposition walk.

Together with estimates on the least eigenvalue, this approach yields $L^2$ mixing time bounds.

To begin, consider a symmetric probability measure $q$ on a finite group $G$ and fix a symmetric set $S$ that generates $G$ and such that $q(s) > 0$ for $s \in S$. Define paths in the Cayley graph $(G, S)$ to be sequences $\delta = (e, y_1, y_2, \ldots, y_k)$, where $e$ is the group identity and $y_i^{-1} y_{i+1} \in S$. Given such a path, define its length to be $|\delta| = k$, and for each $s \in S$, let

$$N(s, \delta) = |\{i \in \{0, \ldots, k-1\} | y_i^{-1} y_{i+1} = s\}|.$$

That is, $N(s, \delta)$ is the number of times the generator $s$ is used in the path $\delta$. Furthermore, let $d_S(x, y)$ denote the graph distance on $(G, S)$ between $x$ and $y$.

DEFINITION 4.1. Fix two symmetric probability measures $\tilde{q}, q$ on a finite group $G$ and a symmetric set generating $S \subset \text{supp}(q)$. A $(\tilde{q}, q)$-flow is a nonnegative function $\eta$ on the set of all paths $\mathcal{P}$ in the Cayley graph $(G, S)$ such that

$$\sum_{\delta \in \mathcal{P}_y} \eta(\delta) = \tilde{q}(y),$$

where $\mathcal{P}_y$ is the set of all paths from the group identity $e$ to $y$ contained in $\mathcal{P}$.

4.1. *The least eigenvalue.* This section presents a lower bound on the smallest eigenvalue of the chain $\frac{1}{2}(q_{n,k_n} + q_{n,k_n}^*)$. The proof relies on a geometric result that bounds the eigenvalues of symmetric chains by considering loops at the identity of odd length. (See [6] for details.) Together with comparison, Lemma 4.1 will be used to derive estimates on mixing time in Section 4.2.

The following definition of an odd flow is analogous to that of a flow, but restricted to paths of odd length.

DEFINITION 4.2. Fix two symmetric probability measures $\tilde{q}, q$ on a finite group $G$ and a symmetric set $S \subset \text{supp}(q)$. A $(\tilde{q}, q)$-odd flow is a nonnegative function $\eta$ on the set of paths of odd length $\mathcal{O}$ in the Cayley graph $(G, S)$ such that

$$\sum_{\delta \in \mathcal{O}_y} \eta(\delta) = \tilde{q}(y),$$

where $\mathcal{O}_y$ is the set of all paths of odd length from the group identity $e$ to $y$ contained in $\mathcal{O}$.

Note that we are not assuming that $S$ generates $G$, that is, the Cayley graph $(G, S)$ need not be connected. However, the existence of a $(\tilde{q}, q)$-odd flow implies that, for each $y$ with $\tilde{q}(y) > 0$, there is at least one path from $e$ to $y$ in $\mathcal{O}$.

THEOREM 4.1 ([6]). *Fix two symmetric probability measures $\tilde{q}, q$ on a group $G$ and a symmetric set $S \subset \text{supp}(q)$. For any $(\tilde{q}, q)$-odd flow $\eta$,*

$$\beta_{\min} \geq -1 + \frac{1 + \tilde{\beta}_{\min}}{A(\eta)},$$

*where $\beta_{\min}$ and $\tilde{\beta}_{\min}$ are the smallest eigenvalues of $q$ and $\tilde{q}$ respectively, and*

$$A(\eta) = \max_{s \in S} \frac{1}{q(s)} \sum_{\delta \in \mathcal{O}} |\delta| N(s, \delta) \eta(\delta).$$

It is well known that a chain $q$ is aperiodic if and only if the least eigenvalue satisfies $\beta_{\min} = -1$. As a trivial application of Theorem 4.1, by taking $S = \{e\}$ and $\tilde{q}(e) = 1$, we have $\beta_{\min} \geq -1 + 2q(e)$. When our chain puts no weight on the identity, the above result provides a way to capture more subtle effects of aperiodicity on the least eigenvalue.

LEMMA 4.1. *Let $\beta_{\min}$ be the smallest eigenvalue of the symmetric chain $\tilde{q}_{n,k_n} = \frac{1}{2}(q_{n,k_n} + q^*_{n,k_n})$. Then*

$$\beta_{\min} \geq -1 + \frac{k_n - 1}{k_n(n - k_n + 2)(n + 1)}.$$

PROOF. We will apply Theorem 4.1 with $\tilde{q}(e) = 1$ and $\tilde{q}(g) = 0$ otherwise. In this case, $\tilde{\beta}_{\min} = 1$. Let $S = \text{supp}(\tilde{q}_{n,k_n})$. For $l$ odd and such that $n - k_n + 1 \leq l \leq n$, define paths

$$\delta_l^{\pm 1} = (e, \sigma_l, \sigma_l^2, \ldots, \sigma_l^l)^{\pm 1}$$

and set $\mathcal{O} = \{\delta_l^{\pm 1} | l \text{ odd}, n - k_n + 1 \leq l \leq n\}$. Let

$$\eta(\delta_l^{\pm 1}) \equiv \frac{1}{2 \sum_{\substack{n-k_n+1 \leq m \leq n \\ m \text{ odd}}} 1/m^2} \cdot \frac{1}{|\delta_l^{\pm 1}|^2}$$

$$\leq \int_{n-k_n+2}^{n+1} \frac{1}{x^2} \cdot \frac{1}{l^2}$$

$$= \frac{(n - k_n + 2)(n + 1)}{k_n - 1} \cdot \frac{1}{l^2}$$

and $\eta(\delta) = 0$ otherwise. Then,

$$A(\eta) \leq \frac{2k_n(n-k_n+2)(n+1)}{k_n-1} \max_{s \in S} \sum_{\delta \in \mathcal{O}} \frac{N(s,\delta)}{|\delta|}$$

$$= \frac{2k_n(n-k_n+2)(n+1)}{k_n-1}.$$

The result follows from Theorem 4.1. $\square$

Theorem 4.1 gives the best results when we can use short paths. In the case of $\frac{1}{2}(q_{n,k_n} + q^*_{n,k_n})$, for paths $\delta$ with $|\delta| \leq \lfloor \frac{n-k_n}{2} \rfloor$, the card originally in position $\lfloor \frac{n-k_n}{2} \rfloor + 1$ moves distance $\pm 1$ at each step along the path. Consequently, the shortest loops at the identity with odd length have length $\approx n - k_n$.

4.2. *Bounds on mixing times.* The following lemma gives a lower bound on the mixing time of $\frac{1}{2}(q_{n,k_n} + q^*_{n,k_n})$ for $k_n$ sufficiently small by looking at the motion of an individual particle.

LEMMA 4.2. *Let $\tilde{q}_{n,k_n} = \frac{1}{2}(q_{n,k_n} + q^*_{n,k_n})$ with $k_n \leq cn$, $0 < c < 1$. Then there is a constant $N(c)$ such that, for $n \geq N$, and $l \leq \frac{c(1-c)^2 n^2}{12}$,*

$$\|\tilde{q}^l_{n,k_n} - \pi\|_{\mathrm{TV}} \geq \frac{c}{2}.$$

*In particular, there is a constant $A(c)$ such that the total variation mixing time satisfies*

$$T(S_n, \tilde{q}_{n,k_n}) \geq An^2.$$

PROOF. Note that the card originally in position $\lfloor \frac{(1-c)n}{2} \rfloor + 1$ performs a simple random walk on $\{1, \ldots, \lfloor (1-c)n \rfloor\}$ before hitting any of the bottom $\lfloor cn \rfloor$ positions. Call this card $a$ and define the event

$$A = \{\sigma | n - \lfloor cn \rfloor < \sigma^{-1}(a) \leq n\},$$

that is, $a$ is in the bottom $\lfloor cn \rfloor$ positions. Then $\pi(A) \geq c - 1/n$. For $l = \lfloor \frac{c(1-c)^2 n^2}{12} \rfloor$, let $X_1, \ldots, X_l$ be an i.i.d. random variable with $P(X_i = \pm 1) = \frac{1}{2}$, and let $S_j = \sum_1^j X_i$. Then

$$\tilde{q}^l_{n,k_n}(A) \leq P\left[\max_{1 \leq j \leq l} |S_j| \geq \frac{(1-c)n}{2}\right]$$

$$\leq \frac{4l}{(1-c)^2 n^2} \qquad \text{(by Kolmogorov's maximal inequality)}$$

$$\leq \frac{c}{3}.$$

Since $\|\tilde{q}^l_{n,k_n} - \pi\|_{\text{TV}} = \max_{A \subset S_n} |\tilde{q}^l_{n,k_n}(A) - \pi(A)|$, the result follows by taking $n$ sufficiently large. The mixing time bound follows from the fact that, for $n$ sufficiently large,

$$c \leq 2\|\tilde{q}^l_{n,k_n} - \pi\|_{\text{TV}} \leq e^{-\lfloor l/T(S_n, \tilde{q}_{n,k_n})\rfloor}.$$

In particular,

$$T(S_n, \tilde{q}_{n,k_n}) \geq \frac{l}{1 - \log c}. \qquad \square$$

Now we will derive an upper bound on the mixing time of $\frac{1}{2}(q_{n,k_n} + q^*_{n,k_n})$ with $n - k_n \leq C$ independent of $n$. That is, the symmetric version of the walk that moves the top card uniformly at random to any but a finite number of the top positions. The proof is by comparison and is based on the following two results. For proofs of these results, see, for example, [4, 5, 6].

DEFINITION 4.3. Given a finite group $G$ and a symmetric probability measure $q$, define the Dirichlet form

$$\mathcal{E}_q(f, f) = \frac{1}{2|G|} \sum_{x,y} |f(xy) - f(x)|^2 q(y).$$

Note that $\mathcal{E}_q(f, f) = ((I - Q)f, f)_{L^2(\pi)}$, where $Qf \mapsto f \star q$ is the Markov operator associated to $q$.

PROPOSITION 4.1. *Assume that $\tilde{\mathcal{E}} \leq A\mathcal{E}$. Then,*

$$T(G, q) \leq T_2(G, q) \preceq \max\left\{AT_2(G, \tilde{q}), A\log|G|, \frac{1}{-\log\beta_-}\right\},$$

*where $\beta_- = \max\{0, -\beta_{\min}\}$.*

THEOREM 4.2. *For any $(\tilde{q}, q)$-flow, $\tilde{\mathcal{E}} \leq A(\eta)\mathcal{E}$ with*

$$A(\eta) = \max_{s \in S} \frac{1}{q(s)} \sum_{\delta \in \mathcal{P}} |\delta| N(s, \delta) \eta(\delta).$$

The proofs of the following two mixing time bounds are by comparison with the random transposition measure on $S_n$:

$$q_{\text{RT},n}(g) = \begin{cases} 1/n, & \text{if } g = e, \\ 2/n^2, & \text{if } g = (i,j), i \neq j, \\ 0, & \text{otherwise.} \end{cases}$$

LEMMA 4.3. *For $k_n \geq n - C$, there exist constants $B(C)$ such that*

$$T_2(S_n, \tilde{q}_{n,k_n}) \leq Bn\log n.$$

PROOF. Let $S = \{\sigma_l^{\pm 1} : n - k_n + 1 \leq l \leq n\}$. First we define paths $\delta_{i,j}$, $1 \leq i < j \leq n$ from $e$ to $(i,j)$ in the Cayley graph $(S_n, S)$:

$$\delta_{i,j} = \begin{cases} \sigma_i^{-1} \sigma_j \sigma_{j-1}^{-1} \sigma_i, & C+1 \leq i < j \leq n, \\ (\sigma_n^{-1})^{C-i+1} \sigma_{C+1}^{-1} \sigma_{j+C-i+1} \sigma_{j+C-i}^{-1} \sigma_{C+1} \sigma_n^{C-i+1}, & \begin{array}{l} 1 \leq i \leq C, \\ i < j \leq n - C, \end{array} \\ (\sigma_{n-C}^{-1})^{C-i+1} \sigma_{C+1}^{-1} \sigma_j \sigma_{j-1}^{-1} \sigma_{C+1} \sigma_{n-C}^{C-i+1}, & \begin{array}{l} 1 \leq i \leq C, \\ j > n - C. \end{array} \end{cases}$$

Define a $(q_{\mathrm{RT},n}, \tilde{q}_{n,k_n})$ flow by $\eta(\delta_{i,j}) = \frac{1}{n^2}$. For $i \leq C$, $|\delta_{i,j}| \leq 2(C-2)$. And each $s \in S$ is used in at most $n$ paths $\delta_{i,j}$ with $i > C$. Consequently,

$$A(\eta) \leq 8[C(C+2)^2 + 1].$$

Since $T(S_n, q_{\mathrm{RT}}) \sim \frac{n}{2} \log n$ (see [7] for details), the result follows by applying Lemma 4.1 and Proposition 4.1 together with Theorem 4.2. $\square$

The following lemma bounds the mixing time of $\tilde{q}_{n,k_n} = \frac{1}{2}(q_{n,k_n} + q_{n,k_n}^*)$ for arbitrary $k_n$. The proof is by comparison with the random transposition measure, but while the flow defined in Lemma 4.3 used only one path for each transposition, here, for most transpositions, we define $k - 1$ paths.

LEMMA 4.4. *There exists a constant $A$ such that the mixing time satisfies*

$$T_2(S_n, \tilde{q}_{n,k_n}) \leq An^3 \log n.$$

PROOF. Let $S = \{\sigma_l^{\pm 1} : n - k_n + 1 \leq l \leq n\}$. The proof is by comparison with the random transposition measure $q_{\mathrm{RT}}$. For $j > i > n - k_n$, define the path

$$\delta_{i,j} = \sigma_i^{-1} \sigma_j \sigma_{j-1}^{-1} \sigma_i.$$

For $i < n - k_n$, we define $k_n - 1$ distinct paths $\delta_{i,j}^l$ with $n - k_n < l < n$. For $j > l$, let

$$\delta_{i,j}^l \equiv (\sigma_l^{-1})^{l-i} \delta_{l,j} \sigma_l^{l-i}$$

and for $i < j \leq l$, let

$$\delta_{i,j}^l \equiv (\sigma_l^{-1})^{l-j} \delta_{l,l+1} (\sigma_l)^{j-i} \delta_{l,l+1} \sigma_l^{j-i} \delta_{l,l+1} \sigma_l^{l-j}.$$

So $|\delta_{i,j}^l| \leq 2n + 12 \leq 3n$. Define a $(q_{\mathrm{RT}}, \tilde{q}_{n,k_n})$-flow by $\eta(\delta_{i,j}) = \frac{1}{n^2}$ and $\eta(\delta_{i,j}^l) = \frac{1}{(k_n-1)n^2}$. Then,

$$A(\eta) \leq \frac{6}{n} \max_s \sum_{\delta_{i,j}^l} N(s, \delta_{i,j}^l) + \frac{8k_n}{n^2} \max_s \sum_{\delta_{i,j}} N(s, \delta_{i,j})$$

$$\leq 18n^2 + \frac{8k_n^2}{n^2}.$$

Since $T(S_n, q_{\text{RT}}) \sim \frac{n}{2} \log n$ (see [7] for details), the result follows by applying Lemma 4.1, Proposition 4.1 and Theorem 4.2. □

The following lemma shows the difficulty in applying path comparison via Theorem 4.2 to bound mixing time.

LEMMA 4.5. *Consider a $(\tilde{q}, q)$-flow $\eta$ on $(G, S)$. For*
$$A(\eta) = \max_{s \in S} \frac{1}{q(s)} \sum_{\delta \in P} |\delta| N(s, \delta) \eta(\delta),$$
*we have the lower bound*
$$A(\eta) \geq \sum_{g \in G} d_S^2(e, g) \tilde{q}(g).$$
*In particular, for $X \subset G$, $A(\eta) \geq d_S^2(e, X) \tilde{q}(X)$.*

PROOF. By averaging over $s$,
$$A(\eta) \geq \sum_{s, \delta} |\delta| N(s, \delta) \eta(\delta)$$
$$= \sum_{\delta} |\delta|^2 \eta(\delta)$$
$$\geq \sum_{g} d_S^2(e, g) \tilde{q}(g).$$
□

Observe that we can always choose a $(\tilde{q}, q)$-flow $\eta$ such that
$$A(\eta) \leq \left( \max_{s \in S} \frac{1}{q(s)} \right) \sum_{g \in G} d_S^2(e, g) \tilde{q}(g).$$

Lemma 4.5 shows that the upper bounds on mixing time that we derive in this section are the best one can do using comparison with the random transposition walk.

Consider a symmetrized variant of the Rudvalis shuffle driven by the measure $r_n$ which is uniform on the generating set $\{\sigma_n, \sigma_n^{-1}, (1, n), id\}$. This walk was analyzed in [15] and an $O(n^3 \log n)$ lower bound was derived for the total variation mixing time (see, e.g., [12] for a matching upper bound). Here we use comparison to extend this result to lower bounds for symmetrized top to bottom-$k$ walks.

LEMMA 4.6. *For $0 < c < 1$ and $k_n \leq cn$, there exists a constant $C > 0$ such that the $L^2$ mixing time satisfies*
$$T_2(S_n, \tilde{q}_{n, k_n}) \geq \frac{Cn^3}{k_n^2} \log n.$$

PROOF. Let $S = \{\sigma_n^{\pm 1}, \tau\}$, where $\tau = (1, n)$ and observe that

$$\sigma_l = \sigma_n \cdot (\sigma_{n-1}^{-1})^{n-l} \cdot \sigma_n^{n-l}$$

$$= \sigma_n \cdot (\sigma_n^{-1} \cdot \tau)^{n-l} \cdot \sigma_n^{n-l}.$$

For $n - k_n < l \leq n$, define paths $\delta_{\sigma_l^{\pm 1}}$ in the Cayley graph $(S_n, S)$ as above, and a corresponding simple $(\tilde{q}_{n,k_n}, r_n)$-flow $\eta$. Then

$$A(\eta) \leq \frac{4}{k_n} \sum_{n-k_n < l \leq n} |\delta_{\sigma_l}|^2$$

$$= \frac{4}{k_n} \sum_{n-k_n < l \leq n} [3(n-l) + 1]^2$$

$$\leq B k_n^2$$

for some universal constant $B$. By Theorem 4.2, $\mathcal{E}_{\tilde{q}_{n,k_n}} \leq B k_n^2 \mathcal{E}_{r_n}$. By Proposition 4.1, together with the lower bound on the mixing time for $\dot{q}_n$ given in [15], we have

$$n^3 \log n \preceq \max\left\{ AT_2(G, \tilde{q}), A \log |G|, \frac{1}{-\log \beta_-} \right\}.$$

By Lemma 4.1, $-1/\log \beta_- = O(n^2)$, and so either $AT_2(G, \tilde{q})$ or $A \log |G|$ is bounded below by $n^3 \log n$. By Lemma 4.2, $n^2 \preceq T_2(G, \tilde{q})$, and so $AT_2(G, \tilde{q}) > A \log |G|$. Consequently, for $n$ sufficiently large,

$$n^3 \log n \leq AT_2(G, \tilde{q})$$

and the result follows. $\square$

**5. Proof of Theorem 1.3.** We show that our estimates on the mixing times for $\tilde{q}_{n,k_n}$ and $q_{n,k_n}$ yield bounds for the lazy top to bottom-$k$ shuffles. In order to transfer mixing time results for the reversible walk $\tilde{q}_{n,k_n}$ to the present case of

$$\hat{q}_{n,k_n} = \tfrac{1}{2}(q_{n,k_n} + \delta_e),$$

we recall the following result.

PROPOSITION 5.1 ([6]). *Let $q$ be a probability measure on $G$ and set $q_* = q \star q^*$. Then*

$$T(G, q) \leq T_2(G, q) \leq 2T_2(G, q_*).$$

*More generally, if $q_v = q^v \star q^{*v}$, then $T_2(G, q) \leq 2v T_2(G, q_v)$. Finally, $q^{*v} \star q^v$ can be used instead of $q_v$.*

LEMMA 5.1. *For $k_n \geq n - C$, there exist constants $B(C)$ such that*
$$T_2(S_n, \hat{q}_{n,k_n}) \leq Bn \log n.$$
*For arbitrary $k_n$, there is a constant $A$ such that*
$$T_2(S_n, \hat{q}_{n,k_n}) \leq An^3 \log n.$$

PROOF. By Proposition 5.1, it is sufficient to prove the bounds for
$$p_{n,k_n} = \hat{q}^*_{n,k_n} \star \hat{q}_{n,k_n}.$$
Observe that
$$\begin{aligned} p_{n,k_n} &= \tfrac{1}{2}(q^*_{n,k_n} + \delta_e) \star \tfrac{1}{2}(q_{n,k_n} + \delta_e) \\ &= \tfrac{1}{2}[\tilde{q}_{n,k_n} + \tfrac{1}{2}(q^*_{n,k_n} \star q_{n,k_n} + \delta_e)] \\ &\geq \tfrac{1}{2}\tilde{q}_{n,k_n}. \end{aligned}$$
Consequently, $\mathcal{E}_{\tilde{q}_{n,k_n}}(f,f) \leq 2\mathcal{E}_{p_{n,k_n}}(f,f)$. Note that $p_{n,k_n}$ is a positive operator and, consequently, has nonnegative eigenvalues. The result then follows from Proposition 4.1, together with the $L^2$ mixing time bounds for $\tilde{q}_{n,k_n}$ derived in Section 4.2. □

To transfer total variation mixing time results for $q_{n,k_n}$ to the lazy top to bottom-$k$ shuffle, we make the following elementary observation.

DEFINITION 5.1. Let $q$ drive a walk on $G$. Then for $p \in (0,1)$, the associated $p$-lazy walk is driven by measure
$$\hat{q}_p = pq + (1-p)\delta_e.$$

LEMMA 5.2. *Let $q$ drive a walk on $G$ with stationary distribution $\pi$, and fix $p, \varepsilon \in (0,1)$. Then there exists a constant $C(p, \varepsilon)$ such that mixing times for $q$ and the associated $p$-lazy walk $\hat{q}_p$ satisfy*
$$T(G, \hat{q}_p) \leq \max\left[\frac{2+\varepsilon}{p}T(G,q), C\right].$$
*Specifically, we can take $C = 80/(p\varepsilon^2)$.*

PROOF. Let $S_m$ be a binomial$(m, p)$ random variable. Then
$$\begin{aligned} \|\hat{q}_p^m - \pi\|_{\mathrm{TV}} &= \frac{1}{2}\sum_{g \in G}|\hat{q}_p^m(g) - \pi(g)| \\ &= \frac{1}{2}\sum_{g \in G}\left|\sum_k P(S_m = k)(q^k(g) - \pi(g))\right| \end{aligned}$$

$$\leq \sum_k P(S_m = k) \cdot \|q^k - \pi\|_{\mathrm{TV}}$$

$$\leq P(S_m \leq 2T(G, q)) + \frac{1}{2e^2}.$$

Taking $\bar{m} \geq \frac{2+\varepsilon}{p} T(G, q)$, by Chebyshev's inequality,

$$P(S_{\bar{m}} \leq 2T(G, q)) \leq P\Big(|S_{\bar{m}} - ES_{\bar{m}}| \geq \Big(1 - \frac{2}{2+\varepsilon}\Big) ES_{\bar{m}}\Big)$$

$$\leq \frac{1-p}{\bar{m}p(1 - 2/(2+\varepsilon))^2}.$$

And consequently,

$$\|\hat{q}_p^{\bar{m}} - \pi\|_{\mathrm{TV}} \leq \frac{1-p}{\bar{m}p(1 - 2/(2+\varepsilon))^2} + \frac{1}{2e^2}$$

$$\leq \frac{1}{2e}$$

for $\bar{m} \geq 80/(p\varepsilon^2)$. □

Now we can transfer the mixing time results for $q_{n,k_n}$ to $\hat{q}_{n,k_n}$.

COROLLARY 5.1. *For $k_n \geq n - \sqrt{(n \log n)/2}$, there exists a constant $C$ such that*

$$T(S_n, \hat{q}_{n,k_n}) \leq Cn \log n.$$

*For $c \in (0, 1)$ and $k_n \geq cn$, there exist constants $A(c)$ such that*

$$T(S_n, \hat{q}_{n,k_n}) \leq An^2 \log^2 n.$$

REMARK 5.1. For $k_n \geq n - \sqrt{(n \log n)/2}$, instead of using Lemma 5.2, we can adapt the coupling of Lemma 3.2 to show $T(S_n, \hat{q}_{n,k_n}) \sim 2n \log n$. The coupling $(X_1^m, X_2^m)$ of $q_{n,k_n}$ yields the coupling

$$(\tilde{X}_1^m, \tilde{X}_2^m) = (X_1^{S_m}, X_2^{S_m})$$

of $\hat{q}_{n,k_n}$, where $S_m$ is an independent binomial$(1/2, m)$ random variable. Then, if $T$ is the coupling time for $(X_1^m, X_2^m)$,

$$P(\tilde{X}_1^m \neq \tilde{X}_2^m) \leq P\Big(S_m \leq \Big(1 + \frac{\varepsilon}{2}\Big) n \log n\Big) + P\Big(T > \Big(1 + \frac{\varepsilon}{2}\Big) n \log n\Big).$$

For $m = 2(1 + \varepsilon)n \log n$, the first term goes to 0 by Chebyshev's inequality, and the second term goes to 0 by the cut-off shown in Lemma 3.2.

The lower bound is also analogous to that given in Lemma 3.2, where we now make the observation that

$$P(\hat{L}_j > m) \geq P\left(L_j > \left(1 - \frac{\varepsilon}{2}\right)n\log n\right) \cdot P\left(S_m \leq \left(1 - \frac{\varepsilon}{2}\right)n\log n\right).$$

So for $k_n \geq n - \sqrt{(n\log n)/2}$ and $\varepsilon \in (0,1)$,

$$\lim_{n \to \infty} \|\hat{q}_{n,k_n}^{(1-\varepsilon)2n\log n} - \pi\|_{\mathrm{TV}} = 1$$

and

$$\lim_{n \to \infty} \|\hat{q}_{n,k_n}^{(1+\varepsilon)2n\log n} - \pi\|_{\mathrm{TV}} = 0.$$

Finally, transferring the lower bounds for $k_n = 2, 3$, which were derived using Wilson's lemma, also requires only a simple argument. Let $\{\eta_i\}$ be i.i.d. Bernoulli random variables with $p = 1/2$, and let $N_t = \sum_{i=1}^{t} \eta_i$. Then if $X_t$ is the top to bottom-$k$ process, the lazy top to bottom-$k$ process is given by $\tilde{X}_t = X_{N_t}$. Using the notation of Theorem 3.2, if $(X_t, Y_t)$ is a lifting of $X_t$, then $(\tilde{X}_t, \tilde{Y}_t) = (X_{N_t}, Y_{N_t})$ is a lifting of $\tilde{X}_t$. It is not hard to check that the assumptions of the theorem are met with $\tilde{\Psi} = \Psi$, $\tilde{\lambda} = 1/2 + 1/2\lambda$, $\tilde{R} = R/2$, and $\tilde{\gamma} = \gamma/2$. Then

$$\frac{\log \tilde{\Psi}_{\max} + (1/2)\log \tilde{\gamma}\varepsilon/(4\tilde{R})}{-\log(1 - \tilde{\gamma})} = \frac{\log \Psi_{\max}/2 + (1/2)\log \gamma\varepsilon/(4R)}{-\log(1 - \gamma/2)}.$$

Using the estimates in Lemma 3.5 and [15], we have the following lower bounds.

COROLLARY 5.2. *For $k_n = 2, 3$ and $\varepsilon > 0$, there exist constants $C(\varepsilon)$, $N > 0$ such that, for $n \geq N$, the lazy top to bottom-$k$ shuffle satisfies*

$$\|\hat{q}_{n,k_n}^m - \pi\|_{\mathrm{TV}} \geq 1 - \varepsilon$$

*for $m \leq Cn^3 \log n$.*

**Acknowledgments.** I thank my Ph.D. adviser Laurent Saloff-Coste for suggesting the topic of this paper, for his helpful advice throughout the project and his continual encouragement.

## REFERENCES

[1] ALDOUS, D. (1983). Random walks on finite groups and rapidly mixing Markov chains. *Séminaire de Probabilités XVII. Lecture Notes in Math.* **986** 243–297. Springer, Berlin. MR770418
[2] DIACONIS, P. (1988). *Group Representations in Probability and Statistics*. IMS, Hayward, CA. MR964069

[3] DIACONIS, P., FILL, J. and PITMAN, J. (1992). Analysis of top to random shuffles. *Combin. Probab. Comput.* **1** 135–155. MR1179244
[4] DIACONIS, P. and SALOFF-COSTE, L. (1993). Comparison techniques for random walk on finite groups. *Ann. Probab.* **21** 2131–2156. MR1245303
[5] DIACONIS, P. and SALOFF-COSTE, L. (1993). Comparison theorems for reversible Markov chains. *Ann. Appl. Probab.* **3** 696–730. MR1233621
[6] DIACONIS, P. and SALOFF-COSTE, L. (1995). Random walks on finite groups: A survey of analytic techniques. In *Probability Measures on Groups and Related Structures* **11** (Z. H. Heyer, ed.) 44–75. World Scientific, River Edge, NJ. MR1414925
[7] DIACONIS, P. and SHAHSHAHANI, M. (1981). Generating a random permutation with random transposition. *Z. Wahrsch. Verw. Gebiete* **57** 159–179. MR626813
[8] DIACONIS, P. and STROOCK, D. (1991). Geometric bounds for eigenvalues for Markov chains. *Ann. Appl. Probab.* **1** 36–61. MR1097463
[9] HILDEBRAND, M. (1990). Rates of convergence of some random processes on finite groups. Ph.D. thesis, Harvard Univ.
[10] JONASSON, J. (2005). Biased random-to-top shuffling. Preprint.
[11] SALOFF-COSTE, L. (1996). Lectures on finite Markov chains. *Lectures on Probability Theory and Statistics. Ecole d'Eté de Probabiltés de Saint-Flour XXVI. Lecture Notes in Math.* **1665** 301–413. Springer, Berlin. MR1490046
[12] SALOFF-COSTE, L. (2004). Random walks on finite groups. In *Probability on Discrete Structures* (H. Kesten, ed.). *Encyclopaedia of Mathematical Sciences* **110** 263–346. Springer, Berlin. MR2023654
[13] SALOFF-COSTE, L. (2004). Total variation lower bounds for finite Markov chains: Wilson's lemma. In *Random Walks and Geometry* 515–532. de Gruyter, Berlin. MR2087800
[14] WILMER, E. (2003). A local limit theorem for a family of non-reversible Markov chains. *J. Theoret. Probab.* **16** 751–770. MR2009201
[15] WILSON, D. B. (2003). Mixing time of the Rudvalis shuffle. *Electron. Comm. Probab.* **8** 77–85. MR1987096
[16] WILSON, D. B. (2004). Mixing times of lozenge tiling and card shuffling Markov chains. *Ann. Appl. Probab.* **14** 274–325. MR2023023

CENTER FOR APPLIED MATHEMATICS
657 RHODES HALL
CORNELL UNIVERSITY
ITHACA, NEW YORK 14853
USA
E-MAIL: sharad@cam.cornell.edu